\documentclass[11pt]{amsart}

\usepackage{amsmath, amscd, amssymb, amsthm, amsfonts}
\usepackage{latexsym}
\usepackage{url}
\usepackage{tikz}
\usetikzlibrary{decorations.markings,decorations.pathreplacing}
\usepackage[T1]{fontenc}
\usepackage{multirow}
\usepackage{comment}
\usepackage{makecell}
\usepackage{colonequals}
\usepackage{calc}
\usepackage{xcolor}
\usepackage{listings}
\usepackage{xcolor}

\tikzset{vertex/.style={circle, draw, inner sep=0pt, minimum size=6pt}}

\numberwithin{equation}{section}

\newtheorem{theorem}{Theorem}[section]
\newtheorem{corollary}[theorem]{Corollary}
\newtheorem{defin}[theorem]{Definition}

\newtheorem{lemma}[theorem]{Lemma}
\newtheorem{proposition}[theorem]{Proposition}
\newtheorem{remark}[theorem]{Remark}
\newtheorem{notation}[theorem]{Notation}

\begin{document}

\title[$\Pi$-shape Trees are $2$-DS] {$\Pi$-shape Trees Are Determined by Their 2-Adjacency Spectra}

\author[Y. Nam]{Yunjae Nam}
\address{Department of Mathematical Sciences\\ KAIST\\ 291 Daehak-ro Yuseong-gu\\ Daejeon 34141\\ South Korea}
\email{namy06354@kaist.ac.kr}

\author[M. Yoon]{Mary Yoon}
\address{Department of Mathematics\\ Korea University\\ 
145 Anam-ro Seongbuk-gu\\ Seoul 02841\\ South Korea}
\email{ymary1201@gmail.com}

\medskip

\begin{abstract}
Determining graphs by their spectra is a central problem in spectral graph theory, yet many tree families remain difficult to characterize using the classical adjacency or Laplacian spectra. 
Recently, the $2$-adjacency matrix was introduced, providing a new spectral invariant.
In this paper, we prove that every $\Pi$-shape graph is determined by its $2$-adjacency spectrum among triangle-free connected graphs, providing a spectral characterization for a graph family that is not determined by the Laplacian spectrum. 
These results demonstrate the effectiveness of the $2$-adjacency spectrum in distinguishing non-isomorphic trees.
\end{abstract}

\subjclass[2020]{05C50,05C05}
\keywords{spectral graph theory, $2$-adjacency matrix, $2$-spectrum, determined by spectrum, $\Pi$-shape graphs}

\maketitle

%%%%%%%%%%%%%%%%%%%%%%%%%%%%%%%%%%%%%%%%%%%%%%%%%%%%%%%%%%%%%%%%%%%%%%%%

\section{Introduction}

In this paper, we let $G$ denote a simple undirected graph with $n$ vertices.
Determining graphs by their spectra is one of the oldest and most fundamental problems in spectral graph theory.
Research on graphs determined by their adjacency, Laplacian, and signless Laplacian spectra has been actively conducted over the past several decades; see \cite{Signless 1, Laplacian 1, Laplacian 2, Signless 2, Laplacian 3, Laplacian 4, Laplacian 5, Van Dam, Signless 3}.
However, these classical spectra often fail to distinguish non-isomorphic graphs, even within restricted families of trees.

Apart from the context of graph determination, the $2$-Laplacian matrix, which can be viewed as a refinement of the Laplacian matrix, was introduced in \cite{Yoon 2-adjacency} alongside the $2$-adjacency matrix.
Let $D_G$ and $A_G$ be the diagonal degree matrix and the adjacency matrix of $G$, respectively.
The $2$-adjacency matrix of $G$, denoted by $A^{(2)}_G$, is defined as
\[
A^{(2)}_G = \frac{1}{12} \left( 16A_G-A_G^2+D_G \right).
\]

We define the $2$-spectrum of $G$ as the multiset consisting of all eigenvalues of $A^{(2)}_G$.
Two non-isomorphic graphs are said to be $2$-cospectral if they have the same $2$-spectrum.
A graph $G$ is said to be determined by its $2$-spectrum (2-DS) if it has no $2$-cospectral mate.
Unlike classical graph matrices that are introduced primarily for combinatorial purposes, both the $2$-adjacency matrix and the $2$-Laplacian matrix were originally motivated by a higher-order approximation of differential operators.
Interestingly, subsequent computational experiments in \cite{Yoon 2-adjacency} demonstrated that this $2$-Laplacian matrix substantially reduces the occurrence of cospectral graphs, showing a significantly lower number of cospectral mates for graphs with at most nine vertices.
Since the $2$-adjacency matrix and the $2$-Laplacian matrix were introduced within the same framework, we also perform computational experiments for the $2$-adjacency matrix using SageMath.
Our computations reveal that every graph with at most nine vertices is determined by its $2$-adjacency spectrum.
The SageMath script used for our verifications is provided in Appendix~\ref{app:sage_code}.
These computational observations suggest that the $2$-adjacency spectrum may provide a more effective invariant than traditional adjacency and Laplacian spectra.
It is therefore natural to investigate whether graph families that are not characterized by classical spectra can be determined by their $2$-adjacency spectra.
To explore this question, we focus on the class of trees with four leaves.

Determining trees by their spectra is famously difficult; a classical result of Schwenk \cite{Schwenk} showed that almost all trees are cospectral with respect to the adjacency matrix.
Thus, even within the class of trees, the adjacency spectrum generally does not provide enough information to distinguish individual trees or tree families.
This makes the determination of specific classes of trees by refined spectral invariants a natural and relevant problem.
We therefore focus first on the class of trees.

A $\Pi$-shape graph is a tree with exactly two vertices of degree $3$, denoted by $\Pi(a,b,c,k,\ell)$, where $b$ represents the distance between the two degree-$3$ vertices, and $a, c, k, \ell$ represent the lengths of the pendant paths attached to them.
Up to isomorphism, we may assume that $a\ge c$, $a\ge k$, $c\ge \ell$, and $k\ge \ell$ if $a=c$.
Note that $a$, $b$, $c$, $k$ and $\ell$ are positive integers. 

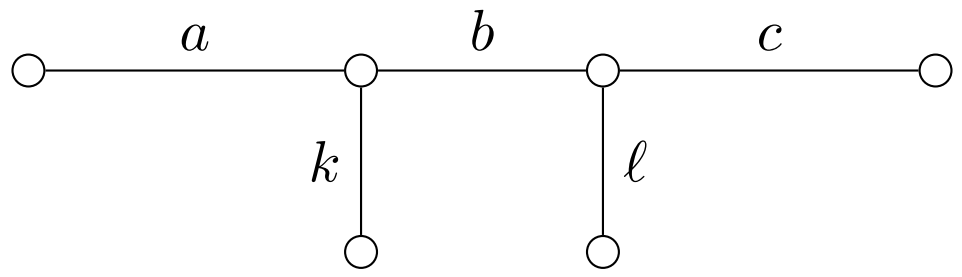
\begin{figure}[htbp]
\centering
\begin{tikzpicture}
\begin{scope}[local bounding box=hshape]
\node[vertex] (L1) at (-3,0) {};
\node[vertex] (L2) at (-0.8,0) {};
\node[vertex] (L3) at (-0.8,-1.2) {};

\node[vertex] (R1) at (3,0) {};
\node[vertex] (R2) at (0.8,0) {};
\node[vertex] (R3) at (0.8,-1.2) {};

\draw (L1) -- (L2) node[midway,above] {$a$};
\draw (L2) -- (L3) node[midway,left] {$k$};

\draw (R1) -- (R2) node[midway,above] {$c$};
\draw (R2) -- (R3) node[midway,right] {$\ell$};

\draw (L2) -- (R2) node[midway,above] {$b$};
\end{scope}

\end{tikzpicture}
\caption{$\Pi$-shape graph}
\label{fig:trees}
\end{figure}

In \cite{pi-shape adjacency}, Liu and Huang classified $\Pi$-shape trees into six types and showed that the family
\[
\Pi(2,1,2,2,2)\cup \Pi(a,b,2,2,2),
\qquad a>2,\ b\ge2,
\]
is determined by its adjacency spectrum.
Furthermore, Wen and Huang proved in \cite{pi-shape Laplacian} that graphs of the form
\[
\Pi(a,1,c,k,\ell)
\]
are mostly determined by their Laplacian spectra, with only a few exceptions, while many examples with $b\ge2$ are not.

Motivated by these results, we investigate whether the finer structural information encoded by the $2$-adjacency spectrum can distinguish $\Pi$-shape trees that are not determined by the classical spectra.
More precisely, we prove that every $\Pi$-shape tree is determined by its $2$-adjacency spectrum among $C_3$-free connected graphs.
In particular, although our result concerns a specific family of trees, the comparison is made with the larger class of all $C_3$-free connected graphs rather than being restricted to trees.
Since the $2$-adjacency matrix of a disconnected graph is block diagonal, its $2$-spectrum is the multiset union of the $2$-spectra of its connected components. Hence, it suffices to consider connected graphs, and we restrict our attention to this case throughout the paper.
Thus, our result shows that the $2$-adjacency spectrum can distinguish $\Pi$-shape trees that are not determined by classical spectra, even when the comparison is made within the broader class of $C_3$-free connected graphs.

The remainder of this paper is organized as follows.
In Section~\ref{sec:comb interp of 2-adj}, we establish a combinatorial interpretation of the traces of powers of the $2$-adjacency matrix and relate them to structural subgraph counts.
In Section~\ref{sec:subgraph counts}, we show that, among $C_3$-free graphs, the $2$-spectrum determines whether a graph is a $\Pi$-shape graph, and derive formulas for the numbers of subgraphs of each possible type in a $\Pi$-shape graph.
Finally, in Section~\ref{sec:determination of pi-shape}, we prove that every $\Pi$-shape graph is determined by its $2$-spectrum among the $C_3$-free connected graphs, thereby establishing the main theorem of this paper.

\section{Combinatorial interpretation of the $2$-Adjacency matrix} \label{sec:comb interp of 2-adj}

Since two graphs have the same $2$-spectrum if and only if the traces of all powers of their $2$-adjacency matrices coincide, we investigate a combinatorial interpretation of the traces of powers of the $2$-adjacency matrix in this section.

\subsection{Definition of the $2$-adjacency matrix}

Recall the definition of $2$-adjacency matrix of $G$.
The $2$-adjacency matrix of $G$, denoted by $A^{(2)}_G$, is defined as
 \[
A^{(2)}_G = \frac{1}{12}\left(16A_{G}-A_{G}^2 + D_{G}\right).
\]

Note that each entry of the matrix $A_G^2 - D_G$ represents the number of open paths of length $2$ between the corresponding vertices of the graph $G$.
More precisely, for $k \ge 1$, let $P_{G,k}$ be the $n \times n$ matrix whose $(i,j)$-entry is equal to the number of open paths of length $k$ from $v_i$ to $v_j$ in $G$. 
Then we can rewrite $A^{(2)}_G$ as follows.

\begin{defin} \label{def-2-adjacency}
    The $2$-adjacency matrix of a simple graph $G$ with $n$ vertices, denoted by $A^{(2)}_G$, is defined as
    \[ 
    A^{(2)}_G = \frac{1}{12}\left(16P_{G,1} - P_{G,2}\right).
    \]
\end{defin}

\medskip

\subsection{Combinatorial interpretation of the trace of $k$th power of $2$-adjacency matrix}

The use of traces of powers of an adjacency matrix to count subgraphs is a well-established approach in spectral graph theory. 
However, in the case of the 2-adjacency matrix, the definition of $A_G^{(2)}$ introduces additional terms involving both $P_{G,1}$ and $P_{G,2}$. 
This reflects the fact that, unlike the ordinary adjacency matrix which encodes only interactions between vertices at distance one, 
the 2-adjacency matrix also incorporates relationships between vertices at distance two. 

As a result, the contributions of $P_{G,1}$ and $P_{G,2}$ appear with different coefficients, 
allowing one to separate and more finely resolve subgraph counts that are indistinguishable in the ordinary adjacency framework. 
This leads to a richer collection of subgraph-counting identities than those obtained from the ordinary adjacency matrix.
In particular, by expressing the traces of powers of $A_G^{(2)}$ in terms of traces of products of $P_{G,1}$ and $P_{G,2}$, 
we obtain relations that allow us to recover structural information about $G$ from its 2-spectrum.

For the remainder of the paper, we focus on the case where $G$ is a tree.
We first introduce some notation. 
For a graph $K$, let $N_G(K)$ denote the number of subgraphs of $G$ that are isomorphic to $K$. 
When no confusion arises, we write $N(K)$ for simplicity.
For each positive integer $m$, let $P_m$ denote the path graph with $m$ vertices.
We denote a starlike tree by $S(a_1,a_2,\ldots,a_m)$, where each $a_i$ is the length of a pendant path attached to a central vertex. 
Without loss of generality, we assume that $a_1 \geq a_2 \geq \cdots \geq a_m\geq 1$.
In particular, we write $S_m=S(1,1,\ldots,1)$ for the star with $m$ leaves.
A $T$-shape tree is a starlike tree with three pendant paths and is written as $S(a,b,c)$.
Moreover, when the graph is of the form $S(a,1,1)$, we denote this graph by $T_{a+2}$.

We now provide a direct combinatorial interpretation of the trace of $(A^{(2)}_G)^k$.
By definition,
\[
\mathrm{tr}\big((A_G^{(2)})^k\big)
=
\sum_{(i_1,\dots,i_k)\in\{1,2\}^k}
\frac{(-1)^{|\{j:i_j=2\}|}16^{|\{j:i_j=1\}|}}{12^k}
\mathrm{tr}(P_{G,i_1}\cdots P_{G,i_k}).
\]
For each $(i_1,\dots,i_k)\in\{1,2\}^k$, let
$W=P_{G,i_1}\cdots P_{G,i_k}$.
Then $W$ is a word of length $k$ in the alphabet $\{P_{G,1},P_{G,2}\}$, and $\mathrm{tr}(W)$ counts the closed walks that follow $W$, with the order of its letters taken into account.

Since
\[
\mathrm{tr}(AB)=\mathrm{tr}(BA),
\]
$\mathrm{tr}(W)$ is invariant under cyclic permutations of the letters of $W$.
Thus, we consider the orbit $W_o=\mathbb{Z}_k\cdot W$ under the cyclic action of
$\mathbb{Z}_k$ on the set of words of length $k$. Let $c_W=|W_o|$ denote the
size of this orbit, and define
\[
\mathrm{tr}(W_o):=\mathrm{tr}(W).
\]
Then
\[
\mathrm{tr}\big((A_G^{(2)})^k\big)
=
\sum_{\mathbb{Z}_k\cdot W}
\frac{(-1)^{|\{j:i_j=2\}|}16^{|\{j:i_j=1\}|}}{12^k}
c_W\mathrm{tr}(W_o).
\]

A walk $w$ is said to follow a word $W$ if, for each
letter $P_{G,i}$ in $W$, the walk makes a move of $i$ steps.
Since the diagonal entries of $P_{G,2}$ are zero, a move corresponding to $P_{G,2}$ cannot start and end at the same vertex.
For a closed walk $w$ following $W$, let $G(w)$ be the induced subgraph of $G$ on the set of vertices visited by $w$. 
We may then decompose $\mathrm{tr}(W_o)$ according to the isomorphism type $K$ of $G(w)$.
Let $n_K(W)$ denote the number of closed walks on $K$ that follow $W$ and generate $K$. 
Since $c_W\mid k$, define $p_W=k/c_W$.
Then
\[
\begin{aligned}
\mathrm{tr}\big((A_G^{(2)})^k\big)
&=
\sum_{\mathbb{Z}_k\cdot W}
\frac{(-1)^{|\{j:i_j=2\}|}16^{|\{j:i_j=1\}|}c_W}{12^k}
\sum_K N_G(K)n_K(W)\\
&=
\frac{k}{12^k}
\sum_{\mathbb{Z}_k\cdot W}\sum_K
\frac{(-1)^{|\{j:i_j=2\}|}16^{|\{j:i_j=1\}|}n_K(W)}
{p_W}N_G(K).
\end{aligned}
\]
Thus, computing $\mathrm{tr}((A_G^{(2)})^k)$ reduces to determining $n_K(W)$.

If $K$ has a vertex $v_0$ with $\deg(v_0)=1$,
then every closed walk following a word and generating $K$ must pass through $v_0$ at least once. 
We may therefore count such walks by fixing $v_0$ as the starting vertex. 
If $W$ is non-periodic and a closed walk passes through $v_0$ exactly $r$ times after leaving the starting point, 
then it gives rise to $r$ distinct closed walks starting at $v_0$ by cyclically shifting the starting point. 
Assigning weight $1/r$ to each such walk therefore
gives total weight $1$ for each closed walk. 
If $W$ is periodic, the same argument gives an additional factor of $p_W$. 
Hence, if $n_{K,r}(W_o)$ denotes the number of closed walks starting at $v_0$,
passing through $v_0$ exactly $r$ times after leaving the starting point,
and following a word in $W_o$,
then
\[
n_K(W)=p_W\sum_r\frac{n_{K,r}(W_o)}{r}.
\]

We introduce the \emph{total length} of $W_o$ by
\[
L(W_o)=L_1(W_o)+2L_2(W_o),
\]
where
\[
L_1(W_o)=|\{j:i_j=1\}|
\qquad\text{and}\qquad
L_2(W_o)=|\{j:i_j=2\}|.
\]
If $K$ is a tree and $L(W_o)$ is odd, then $n_K(W)=0$. 
Since every subgraph of a tree is also a tree, we obtain, 
by dividing
$\mathbb{Z}_k\cdot W$ into cases according to $L_1(W_o)=2i$,
\[
\mathrm{tr}\big((A_G^{(2)})^k\big)
=
k\left(-\frac{1}{12}\right)^k
\sum_{i=0}^{\lfloor k/2\rfloor}
16^{2i}
\sum_{K\text{: a tree}}
\left[
\sum_{L_1(W_o)=2i}\sum_r
\frac{n_{K,r}(W_o)}{r}
\right]N_G(K).
\]

Let \[
    c_{2i}^k(K) := \sum_{L_1(W_o)=2i}\sum_r \frac{n_{K,r}(W_o)}{r}.
\]
Then 
\[
\mathrm{tr}\big((A_G^{(2)})^k\big)
=
k\left(-\frac{1}{12}\right)^k
\sum_{i=0}^{\lfloor k/2\rfloor}
16^{2i}
\sum_{K\text{: a tree}}c_{2i}^k(K)N_G(K).
\]
Thus, our problem reduces to determining $c_{2i}^k(K)$ for each $(i,K)$.

We first provide formulas for $c^{k}_{2i}(K)$ when $K=P_m$, $S_m$, or $T_n$ in the following proposition.
For $k=1,2,3,4$, we then compute $c^{k}_{2i}(K)$ for all remaining possible subgraphs $K$ and express
$\mathrm{tr}\big((A_G^{(2)})^k\big)$ in terms of the numbers of subgraphs $N_G(K)$.

\begin{proposition} \label{prop:coefficient}
For $n\geq 0$ and $0\leq j\leq n$, we have
\[
c_{2j}^{n+j}(P_{n+1})
=
\sum_{i=0}^{n-j}
\binom{n-i}{i}
\binom{j+i}{n-j-i}.
\]

Moreover, for $n\geq 2$ and $0\leq j\leq n-1$,
\[
c_{2j}^{n+j}(S_n)
=
\binom{n}{j}(n-1)!.
\]
Finally, for $n\geq 3$ and $j\geq 0$,
\[
c_{2j}^{n+j}(T_{n})
=
2\sum_{k=0}^{1}\sum_{i=0}^{n-j-k}
\binom{n-1-i}{i}
\binom{j+k+i-1}{n-j-k-i},
\]
where we use the convention that $\binom{n}{m}=0$ if $m<0$ or $m>n$.
\end{proposition}

\begin{proof}
When $K$ has $n$ edges and $j=n$, $c_{2n}^{2n}(K)$ is the number of all walks that start at an arbitrary leaf and traverse each edge exactly twice.
In the case of $K = P_{n+1}$, computing $c_{2j}^{n+j}(K)$ reduces to counting the number of ways to travel from one leaf to the other and then return to the starting leaf, using $n-j$ many $2$-steps and $2j$ many $1$-steps.
Therefore, if one side contains $i$ many $2$-steps, the other side must contain $n-j-i$ many $2$-steps. 
Consequently,
\[
c_{2j}^{n+j}(P_{n+1})
=
\sum_{i=0}^{n-j}
\binom{n-i}{i}
\binom{j+i}{n-j-i}.
\]
In this case, if $n$ is odd and $j=0$, then at least one of the binomial coefficients
\[
\binom{n-i}{i} \quad \text{and} \quad \binom{j+i}{n-j-i}
\]
must vanish for every valid $i$, and hence the sum evaluates to $0$.

In the case of $K=S_n$, we have $c_{2n}^{2n}(S_n)=(n-1)!$, and for each of the $(n-1)!$ cases, the closed walk consists of $n$ paths of length $2$.
In this case, we choose $n-j$ of the $n$ paths and replace two $1$-steps by one $2$-step, 
so there are $\binom{n}{n-j}=\binom{n}{j}$ possibilities for each case.

In the case of $K=T_n$, both of the $2$ cases consist of two paths of length $n-1$ and one path of length $2$.
Dividing the possibilities into the case where the path of length $2$ is replaced by one $2$-step and the case where it is left as two $1$-steps, we immediately obtain
\[
c_{2j}^{n+j}(T_{n})
=
2\sum_{k=0}^{1}\sum_{i=0}^{n-j-k}
\binom{n-1-i}{i}
\binom{j+k+i-1}{n-j-k-i}.\qedhere
\]
\end{proof}

Proposition~\ref{prop:coefficient} allows us to determine the smallest power of $A_G^{(2)}$ whose trace contains a contribution from a given subgraph.

\begin{corollary}\label{cor:first occurrence}
For each nonnegative integer $m$, the smallest integer $k$ for which the expansion of $\mathrm{tr}\bigl((A_G^{(2)})^k\bigr)$ contains a contribution from $P_{2m+2}$ and $P_{2m+3}$ is $k=2m+2$. 
Moreover, the coefficients of $N_G(P_{2m+2})$ and $N_G(P_{2m+3})$ in this first occurrence are, respectively,
\[
(2m+2)\left(-\frac{1}{12}\right)^{2m+2}16^2(m+1)^2
\quad\text{and}\quad
(2m+2)\left(-\frac{1}{12}\right)^{2m+2}.
\]
Similarly, for the subgraphs $T_{2m+2}$ and $T_{2m+3}$, the smallest such power is $k=2m+3$.
The corresponding coefficients of $N_G(T_{2m+2})$ and $N_G(T_{2m+3})$ are, respectively,
\[
2(2m+3)\left(-\frac{1}{12}\right)^{2m+3}16^2(m+2)^2
\quad\text{and}\quad
2(2m+3)\left(-\frac{1}{12}\right)^{2m+3}.
\]
Furthermore, for $S(m_1,m_2,m_3)$ with 
\[
m_1+m_2+m_3=m,
\]
the smallest power $k$ yielding a contribution from $S(m_1,m_2,m_3)$ is $k=m$ if $m_1,m_2$, and $m_3$ are either all even or all odd. 
Otherwise, the smallest power is $k=m+1$.
More precisely, the coefficient of $N_G(S(m_1,m_2,m_3))$ is given by
\[
2(m+1)\left(-\frac{1}{12}\right)^{m+1}16^2
\prod_{\substack{r\in\{1,2,3\}\\ m-m_r\text{ is odd}}}
\left\lceil\frac{m-m_r}{2}\right\rceil
\]
if $m-m_r$ is odd for some $r$, and by
\[
2m\left(-\frac{1}{12}\right)^m
\]
if $m-m_r$ is even for all $r$.
\end{corollary}

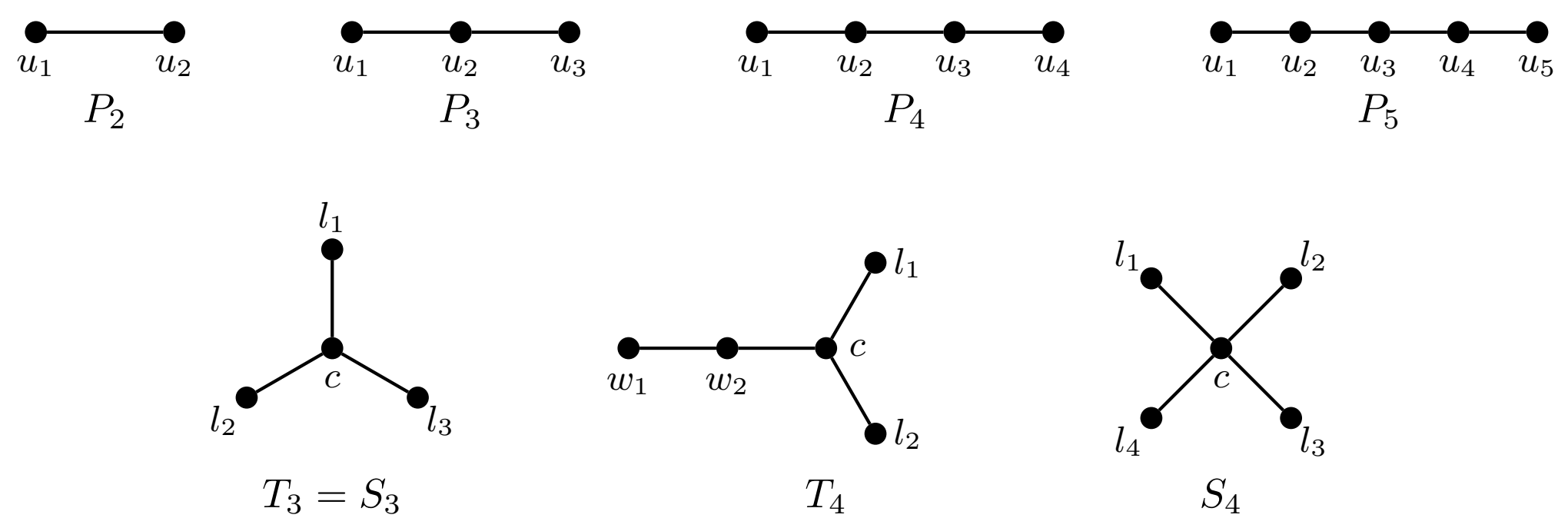
\begin{figure}[htbp]
\centering
\begin{tikzpicture}[scale=0.85, every node/.style={circle, draw, fill=black, inner sep=0pt, minimum size=5pt}]

    % =================================================================
    % ROW 1: Path-type Supporting Subtrees (P_2, P_3, P_4, P_5)
    % =================================================================
    
    % --- P_2 ---
    \begin{scope}[shift={(0,5.2)}]
        \node (u1) at (0,0) [label=below:\footnotesize $u_1$] {};
        \node (u2) at (1.4,0) [label=below:\footnotesize $u_2$] {};
        \draw[thick] (u1) -- (u2);
        \node[draw=none, fill=none] at (0.7,-0.8) {\small $P_2$};
    \end{scope}

    % --- P_3 ---
    \begin{scope}[shift={(3.2,5.2)}]
        \node (u1) at (0,0) [label=below:\footnotesize $u_1$] {};
        \node (u2) at (1.1,0) [label=below:\footnotesize $u_2$] {};
        \node (u3) at (2.2,0) [label=below:\footnotesize $u_3$] {};
        \draw[thick] (u1) -- (u2) -- (u3);
        \node[draw=none, fill=none] at (1.1,-0.8) {\small $P_3$};
    \end{scope}

    % --- P_4 ---
    \begin{scope}[shift={(7.3,5.2)}] 
        \node (u1) at (0,0) [label=below:\footnotesize $u_1$] {};
        \node (u2) at (1.0,0) [label=below:\footnotesize $u_2$] {};
        \node (u3) at (2.0,0) [label=below:\footnotesize $u_3$] {};
        \node (u4) at (3.0,0) [label=below:\footnotesize $u_4$] {};
        \draw[thick] (u1) -- (u2) -- (u3) -- (u4);
        \node[draw=none, fill=none] at (1.5,-0.8) {\small $P_4$};
    \end{scope}

    % --- P_5 ---
    \begin{scope}[shift={(12.0,5.2)}]
        \node (u1) at (0,0) [label=below:\footnotesize $u_1$] {};
        \node (u2) at (0.8,0) [label=below:\footnotesize $u_2$] {};
        \node (u3) at (1.6,0) [label=below:\footnotesize $u_3$] {};
        \node (u4) at (2.4,0) [label=below:\footnotesize $u_4$] {};
        \node (u5) at (3.2,0) [label=below:\footnotesize $u_5$] {};
        \draw[thick] (u1) -- (u2) -- (u3) -- (u4) -- (u5);
        \node[draw=none, fill=none] at (1.6,-0.8) {\small $P_5$};
    \end{scope}

    % =================================================================
    % ROW 2: Branching-type Supporting Subtrees (T_3, S(1,1,1,1))
    % =================================================================

    % --- T_3 (Claw) ---
    \begin{scope}[shift={(3,2)}]
        \node (c) at (0,0) [label={[label distance=2pt]below:\footnotesize $c$}] {};
        \node (l1) at (90:1) [label=above:\footnotesize $l_1$] {};
        \node (l2) at (210:1) [label=below left:\footnotesize $l_2$] {};
        \node (l3) at (-30:1) [label=below right:\footnotesize $l_3$] {};
        \draw[thick] (c) -- (l1);
        \draw[thick] (c) -- (l2);
        \draw[thick] (c) -- (l3);
        \node[draw=none, fill=none] at (0,-1.5) {\small $T_3=S_3$};
    \end{scope}

    % --- S(1,1,1,1) ---
    \begin{scope}[shift={(12,2)}]
        \node (c) at (0,0) [label={[label distance=2pt]below:\footnotesize $c$}] {};
        \node (l1) at (135:1) [label=above left:\footnotesize $l_1$] {};
        \node (l2) at (45:1) [label=above right:\footnotesize $l_2$] {};
        \node (l3) at (-45:1) [label=below right:\footnotesize $l_3$] {};
        \node (l4) at (-135:1) [label=below left:\footnotesize $l_4$] {};
        \draw[thick] (c) -- (l1);
        \draw[thick] (c) -- (l2);
        \draw[thick] (c) -- (l3);
        \draw[thick] (c) -- (l4);
        \node[draw=none, fill=none] at (0,-1.5) {\small $S_4$};
    \end{scope}

    % --- T_4 ---
    \begin{scope}[shift={(8,2)}]
        \node (c) at (0,0) [label={[label distance=2pt]right:\footnotesize $c$}] {};
        \node (l1) at (60:1) [label=right:\footnotesize $l_1$] {};
        \node (l2) at (-60:1) [label=right:\footnotesize $l_2$] {};
        \node (w2) at (-1,0) [label=below:\footnotesize $w_2$] {};
        \node (w1) at (-2,0) [label=below:\footnotesize $w_1$] {};
        \draw[thick] (c) -- (l1);
        \draw[thick] (c) -- (l2);
        \draw[thick] (c) -- (w2);
        \draw[thick] (w1) -- (w2);
        \node[draw=none, fill=none] at (0,-1.5) {\small $T_4$};
    \end{scope}

\end{tikzpicture}
\label{fig:totallengthleq8}
\caption{All supporting subtrees arising in the expansion of $(A_G^{(2)})^4$,
i.e., those generated by closed walks following $W_o$ with $L(W_o)\leq 8$.}
\end{figure}

\begin{proposition} \label{k-power trace}
    Let $G$ be a tree. Then
\begin{align*}
\mathrm{tr}(A_{G}^{(2)}) &= 0,\\
\mathrm{tr}\big((A_{G}^{(2)})^2\big)
&= \frac{2}{12^2}\left(16^2N(P_2)+N(P_3)\right),\\
\mathrm{tr}\big((A_{G}^{(2)})^3\big)
&= \frac{3}{12^3}\left(-2\cdot16^2N(P_3)-2N(T_3)\right),\\
\mathrm{tr}\big((A_{G}^{(2)})^4\big)
&= \frac{4}{12^4}\Big(
    16^4\big(\frac{1}{2}N(P_2)+N(P_3)\big) \\
&\quad + 16^2\big(2N(P_3)+4N(P_4)+6N(T_3)\big) \\
&\quad + \frac{1}{2}N(P_3)+N(P_5)+3N(T_3)+6N(S_4) 
\Big).
\end{align*}
\end{proposition}
\begin{proof}
Since $K$ is a tree, every edge of $K$ must be traversed at least twice by a closed walk generating $K$. 
Hence,
$n_K(W_o)=0$ if $L(W_o)<2E(K)$.
Therefore, when $k\leq 4$, we have $L(W_o)\leq 2k\leq 8$, so it suffices to consider trees with at most $5$ vertices.
The trees with at most $5$ vertices are precisely
\[
    P_2,P_3,P_4,P_5,T_3=S_3,T_4,S_4.
\]
See Figure~\ref{fig:totallengthleq8} for these trees.
When $k\leq 4$, the remaining cases are those for which $L(W_o)>2E(K)$ and hence are not covered above. They are as follows.
\begin{enumerate}
\item For $K=P_2$, the walk must consist entirely of $1$-step moves. Thus, the only remaining case that can yield a nonzero $c_{2i}^j(K)$ is
\[
(L_1(W_o),L_2(W_o))=(4,0).
\]
In this case, the unique closed walk traverses the edge back and forth twice, so
\[
c_4^4(P_2)=\frac{1}{2}.
\]
\item For $K=P_3$, the only cases not covered above that can yield a nonzero $c_i^j(K)$ are
\[
(L_1(W_o),L_2(W_o))=(0,3),(2,2),(0,4),
\]
since these are the cases with $L(W_o)>4$. There is no possible walk in the case $(0,3)$, so it remains to consider the other two cases.\\
In the case $(2,2)$, the only possible walks are obtained by inserting two $1$-step moves into the corresponding walk with two $2$-steps.  
We write
\[ u \xrightarrow{v} w \]
for a $2$-step move from $u$ to $w$ through $v$.
Therefore, %writing \yunjae{Erase here}
%\[
%P_3=u_1-u_2-u_3,
%\]
the possible walks are
\[
\begin{aligned}
&u_1\xrightarrow{u_2} u_3\to u_2\to u_3\xrightarrow{u_2}u_1, \quad (r=1)\\
&u_1\to u_2\to u_1\xrightarrow{u_2} u_3\xrightarrow{u_2}u_1, \quad (r=2)\\
&u_1\xrightarrow{u_2} u_3\xrightarrow{u_2}u_1\to u_2\to u_1. \quad (r=2)
\end{aligned}
\]
Hence,
\[
c_2^4(P_3)=2.
\]
Furthermore, in the case $(0,4)$, there is a unique possible move at each step, and therefore
\[
c_0^4(P_3)=1.
\]
\item For $K=T_3$, the only possible case is
\[
(L_1(W_o),L_2(W_o))=(0,4).
\]
There are two possibilities in which the walk visits only the two leaves other than the starting point:
\[
l_1\xrightarrow{c} l_2\xrightarrow{c}  l_3\xrightarrow{c}  l_2\xrightarrow{c} l_1
\qquad\text{and}\qquad
l_1\xrightarrow{c} l_3\xrightarrow{c}  l_2\xrightarrow{c}  l_3\xrightarrow{c} l_1,
\]
and two possibilities in which the walk visits $l_1$ twice after leaving the starting point, i.e. $r=2$,
\[
l_1\xrightarrow{c} l_2\xrightarrow{c}  l_1\xrightarrow{c}  l_3\xrightarrow{c} l_1
\qquad\text{and}\qquad
l_1\xrightarrow{c} l_3\xrightarrow{c}  l_1\xrightarrow{c}  l_2\xrightarrow{c} l_1,
\]
Hence,
\[
c_0^4(T_3)=3.
\]
\end{enumerate}
\end{proof}

\medskip

\section{Subgraph Counts in $\Pi$-Shape Graphs} \label{sec:subgraph counts}

In this section, we show that, among $C_3$-free graphs, the $2$-spectrum determines whether a graph is a $\Pi$-shape. 
We also derive formulas for the numbers of subgraphs of each possible type in a $\Pi$-shape graph. 
These formulas will be used in the next section to show that $\Pi$-shape graphs are determined by their $2$-spectra.

\begin{lemma} \label{lem:determine H of tree}
Let $G$ be a tree with $n$ vertices. 
Then the $2$-spectrum of $G$ determines whether $G$ is isomorphic to a graph of the form $\Pi(a,b,c,k,\ell)$.
\end{lemma}
\begin{proof}
Since $G$ is a tree, we have $N(P_2)=n-1$.
By Proposition~\ref{k-power trace}, the $2$-spectrum determines both $N(P_3)$ and $N(T_3)$.
Moreover,
\[
N(T_3)=\sum_{v\in V(G)}\binom{d(v)}{3},
\]
where $d(v)$ denotes the degree of $v$.
Thus, $N(T_3)=2$ implies that exactly two vertices of $G$ have degree $3$, while all other vertices have degree at most $2$.
Hence, $G$ is a $\Pi$-shape graph if and only if $N(T_3)=2$.
\end{proof}

\begin{lemma}  \label{lem:determine H of C3 free}
Let $G$ be a $\Pi$-shape graph and $G'$ be a $C_3$-free connected graph. 
If $G$ and $G'$ are $2$-cospectral, then $G'$ is also a $\Pi$-shape graph.
\end{lemma}
\begin{proof}
Since $G$ and $G'$ are $2$-cospectral mates, we have
\[
\mathrm{tr}\big((A_G^{(2)})^2\big)
=
\mathrm{tr}\big((A_{G'}^{(2)})^2\big).
\]
By Proposition~\ref{k-power trace}, 
\[
\mathrm{tr}\big((A_{G'}^{(2)})^2\big)
=
\frac{2}{12^2}
\left(
16^2N_{G'}(P_2)
+N_{G'}(P_3)
+4N_{G'}(C_4)
\right).
\]

Since $G'$ is connected,
\[
N_{G'}(P_2)\geq n-1
\quad\text{and}\quad
N_{G'}(P_3)\geq n-2.
\]
Suppose, for a contradiction, that $G'$ is not a tree. Then
$N_{G'}(P_2)\geq n$, and hence
\[
\mathrm{tr}\big((A_{G'}^{(2)})^2\big)
\geq
\frac{2}{12^2}
\left(
16^2n
+N_{G'}(P_3)
+4N_{G'}(C_4)
\right).
\]

On the other hand, since $G$ is a $\Pi$-shape graph and hence a tree,
Proposition~\ref{k-power trace} gives
\[
\mathrm{tr}\big((A_G^{(2)})^2\big)
=
\frac{2}{12^2}
\left(
16^2(n-1)+n
\right).
\]
Therefore, from
\[
\mathrm{tr}\big((A_G^{(2)})^2\big)
=
\mathrm{tr}\big((A_{G'}^{(2)})^2\big),
\]
we obtain
\[
16^2(n-1)+n
\geq
16^2n
+N_{G'}(P_3)
+4N_{G'}(C_4).
\]
Thus,
\[
N_{G'}(P_3)+4N_{G'}(C_4)
\leq n-256.
\]
Since $N_{G'}(P_3)\geq n-2$, this is impossible. Hence, $G'$ must be a tree.
Thus, by Lemma~\ref{lem:determine H of tree}, $G'$ is a $\Pi$-shape graph.
\end{proof}

For the remainder of this paper, we regard $G$ as a graph of the form
\[
\Pi(a,b,c,k,\ell).
\]

Before stating the following proposition, we introduce some notations. 
Whenever a collection of parameters is written using set notation, we regard it as a multiset, with multiplicities taken into account.

\begin{notation}
For a $\Pi$-shape graph $\Pi(a,b,c,k,\ell)$ and each positive integer $i$,
define
\[
\begin{aligned}
\alpha_i
&=
\left|\{x\in\{a,c,k,\ell\}:x=i\}\right|,\\
\beta_i
&=
\begin{cases}
1,&b=i,\\
0,&\text{otherwise},
\end{cases}\\
\gamma_i
&=
\left|\{x\in\{a+k,c+\ell\}:x=i\}\right|,\\
A_i
&=
\left|\{x\in
\{a+b+c,a+b+\ell,k+b+\ell,k+b+c\}:x=i\}\right|,\\
B_i
&=
\left|\{x\in
\{a+b,b+c,b+k,b+\ell\}:x=i\}\right|.
\end{aligned}
\]
\end{notation}

For a $\Pi$-shape graph $\Pi(a',b',c',k',\ell')$,
we use the notation
\[
\alpha_i',\ \beta_i',\ \gamma_i',\ A_i',\ \text{and } B_i'
\]
for the corresponding quantities.

Every connected subgraph of a $\Pi$-shape graph is of one of the following three types: a path, a $T$-shape, or a $\Pi$-shape.
To describe the numbers of subgraphs of these three types, let $R_m(n)$, $U_m(n)$, and $\delta_t(n)$ denote the ramp function, unit step function, and delta function, respectively:
\[
R_m(n)=
\begin{cases}
n-m, & \text{if } n\geq m,\\
0, & \text{if } n<m,
\end{cases}
\]
\[
U_m(n)=
\begin{cases}
1, & \text{if } n\geq m,\\
0, & \text{if } n<m,
\end{cases}
\]
\[
\delta_t(n)=
\begin{cases}
1, & \text{if } t=n,\\
0, & \text{if } t\neq n.
\end{cases}
\]

Using these functions, we can conveniently express the numbers of subgraphs of the three types.

\begin{lemma} [Number of path subgraphs] \label{lem:num of path} 
    Let $G$ be a graph with $n$ vertices isomorphic to $\Pi(a,b,c,k,\ell)$. 
    Then, for each integer $m \geq 1$, 
\[ N(P_{m+3}) = n + m - \sum_{i=1}^m i\,(\alpha_{m+1-i}-\beta_{m+1-i}-\gamma_{m+1-i}-A_{m+1-i}+B_{m+1-i}).\]
\end{lemma}
\begin{proof}
Using the inclusion--exclusion principle, we obtain
\begin{align*}
N(P_{m+3})
=\,&\ R_{m+1}(a+b+c)+ R_{m+1}(k+b+c)+ R_{m+1}(a+b+\ell)+ R_{m+1}(k+b+\ell)
\\
&+ R_{m+1}(a+k)+ R_{m+1}(c+\ell)
\\
&-R_{m+1}(b+c)-R_{m+1}(a+b)-R_{m+1}(k+b)-R_{m+1}(b+\ell)
\\
&-R_{m+1}(a)-R_{m+1}(c)-R_{m+1}(\ell)-R_{m+1}(k)
\\
&+R_{m+1}(b).
\end{align*}
From this expression, the desired formula follows.
\end{proof}

\begin{lemma}[Number of $T$-shape subgraphs]\label{lem:num_tshape} 
The number of subgraphs of $\Pi(a,b,c,k,\ell)$ isomorphic to
    $S(m_1,m_2,m_3)$ is
    \begin{align*}
        &\sum_{(r_1,r_2,r_3)\in M_3}
        U_{r_1}(a)
        U_{r_2}(k)
        \big(
        U_{r_3}(b+c)
        +U_{r_3}(b+\ell)
        -U_{r_3}(b)
        \big)\\
        &\quad+
        \sum_{(r_1,r_2,r_3)\in M_3}
        U_{r_1}(c)
        U_{r_2}(\ell)
        \big(
        U_{r_3}(a+b)
        +U_{r_3}(k+b)
        -U_{r_3}(b)
        \big).
    \end{align*}
   Here, $M_3$ denotes the set of distinct tuples obtained by all permutations of $(m_1, m_2, m_3)$ (i.e., the set of distinct elements in $\{ (m_{\sigma(1)}, m_{\sigma(2)}, m_{\sigma(3)}) : \sigma \in \mathfrak{S}_3 \}$).
    In particular, when $m_1=m+1$ and $m_2=m_3=1$,
    \[
        N(T_{m+3})
        =6-\sum_{i=1}^m(\alpha_i-2\beta_i+B_i).
    \]
\end{lemma}
\begin{proof}
A $\Pi$-shape graph has exactly two vertices of degree $3$, say $v$ and $w$, 
while all other vertices have degree at most $2$. 
Thus, any $T$-shape subgraph must contain either $v$ or $w$ as its unique vertex of degree $3$.

Suppose first that the degree-$3$ vertex is $v$.
Two of the pendant paths of the $T$-shape must lie along the pendant paths of lengths $a$ and $k$ at $v$.
For a fixed assignment $(r_1,r_2,r_3)\in M_3$ of the three pendant-path lengths, the first two paths can be chosen in
\[
U_{r_1}(a)U_{r_2}(k)
\]
ways.
The third path can either end on the pendant path of length $c$ or on the pendant path of length $\ell$, after passing through the path of length $b$.
The number of possible choices for the third path is therefore
\[
U_{r_3}(b+c)+U_{r_3}(b+\ell)-U_{r_3}(b),
\]
where the last term corrects for the common possibility that the third path ends at $w$.
This gives the first sum.

The case in which the degree-$3$ vertex is $w$ is analogous and gives the second sum.
\end{proof}

\begin{lemma}[Number of $\Pi$-shape subgraphs]\label{lem:num_four_leaves}
Let $G\cong\Pi(a,b,c,k,\ell)$.
Then the number of subgraphs of $G$ isomorphic to
$\Pi(m_1,t,m_2,m_3,m_4)$ is
\[
\delta_t(b)
\sum_{(r_1,r_2,r_3,r_4)\in M_4}
U_{r_1}(a)
U_{r_2}(c)
U_{r_3}(k)
U_{r_4}(\ell).
\]
Here, $M_4$ denotes the set of distinct tuples generated by the action of $\mathfrak{D}_4$ on $(m_1, m_2, m_3, m_4)$, 
where
$\mathfrak{D}_4=
\langle (1\,3),(2\,4),(1\,2)(3\,4)\rangle$.

\end{lemma}
\begin{proof}
A $\Pi$-shape graph has exactly two vertices of degree $3$, say $v$ and $w$, 
while all other vertices have degree at most $2$. 
Hence, any $\Pi$-shape subgraph must contain both $v$ and $w$ as its two degree $3$ vertices. 
Since the distance between $v$ and $w$ in $G$ is $b$, such a subgraph can occur only when $t=b$, 
which accounts for the factor $\delta_t(b)$.

The four pendant paths incident to $v$ and $w$ can be assigned to the four pendant paths of $\Pi(m_1,t,m_2,m_3,m_4)$ in precisely the ways given by $\mathfrak{D}_4$. 
For each $(r_1,r_2,r_3,r_4)\in M_4$, the corresponding path lengths can be chosen independently, 
and the number of choices is given by
\[
U_{r_1}(a)
U_{r_2}(c)
U_{r_3}(k)
U_{r_4}(\ell).
\]
Summing over all $\sigma\in\mathfrak{D}_4$ gives the stated formula.
\end{proof}

\medskip

\section{$2$-Spectral Determination of $\Pi$-Shape Graphs}\label{sec:determination of pi-shape}

In Section~\ref{sec:subgraph counts}, we introduced the quantities $\alpha_i$, $\beta_i$, $\gamma_i$, $A_i$, and $B_i$ for a $\Pi$-shape graph $\Pi(a,b,c,k,\ell)$. 
We also used this notation to derive formulas for the numbers of its subgraphs.

In this section, we examine precisely what information each of $\alpha_i$, $\beta_i$, $\gamma_i$, $A_i$, and $B_i$ encodes about
$\Pi(a,b,c,k,\ell)$. 
In particular, we show that the $2$-spectrum determines the numbers of $P_j$ and $T_j$ subgraphs of $\Pi(a,b,c,k,\ell)$ for every $j$. 
Finally, by combining these results, we prove our main theorem: every $\Pi$-shape graph is determined by its $2$-spectrum among the $C_3$-free graphs.

\begin{remark} \label{rmk:a,c,k,l}
Note that $\alpha_j$ and $\gamma_j$ for $j=1,2,\dots,m$ determine precisely which of the parameters $a,c,k,\ell$ belong to $\{1,2,\dots,\lfloor m/2\rfloor\}$ whenever at least two of them belong to this set. 
If exactly one of them belongs to the set, then $\alpha_j$ and $\gamma_j$ may not distinguish whether it is $k$ or $\ell$. Meanwhile, $\beta_j$ for $j=1,2,\dots,m$ determines whether $b=1$, $b=2$, $\dots$, $b=m$, or $b>m$.
\end{remark}

\begin{lemma} \label{lem:A,B}
Suppose that $\alpha_j$, $\beta_j$, and $\gamma_j$ are known for all $j=1,2,\dots,m$. 
Then $A_{m+1}$ and $B_{m+1}$ are uniquely determined.
\end{lemma}

\begin{proof}
We can determine from the $\beta_j$'s whether $b\le m$ or not. 
If $b\geq m+1$, then $A_{m+1}=0$. 
Thus, assume that $b \le m$. 
Then
\[
A_{m+1}
=
\left|
\left\{
x\in
\{a+c,a+\ell,k+\ell,k+c,a+k,c+\ell\}:x=m+1-b
\right\}
\right|
-\gamma_{m+1-b}.\]
Since the $\alpha_j$'s determine, for each $j=1,\ldots,m$, the number of elements among $a,c,k,\ell$ that are equal to $j$, we can determine $A_{m+1}$.

On the other hand,
\[
B_{m+1}
=
\sum_{i=1}^{m}\beta_{m+1-i}\alpha_i
=
\alpha_{m+1-b}.
\]
Since the index $m+1-b$ is at most $m$, 
$B_{m+1}$ is determined by the given $\alpha_j$'s.
\end{proof}

Using Lemmas~\ref{lem:num_tshape}, \ref{lem:num_four_leaves}, and Remark~\ref{rmk:a,c,k,l}, we obtain the following.
\begin{lemma}\label{lem:Pi and T determined}
Let $G\cong\Pi(a,b,c,k,\ell)$. Then the following hold.
\begin{enumerate}
\item If $\alpha_i$, $\beta_i$, and $\gamma_i$ are known for $i=1,\dots,2m$, then
\[N_G\big(\Pi(m_1,t,m_2,m_3,m_4)\big)\]
is determined for all
$m_1+m_2+m_3+m_4+t\leq 2m+4$.

\item If $\alpha_i$, $\beta_i$, and $\gamma_i$ are known for $i=1,\dots,2m+1$, then
\[
N_G\big(\Pi(m_1,t,m_2,m_3,m_4)\big)
\]
is determined for all
$m_1+m_2+m_3+m_4+t\leq 2m+4$,
and $N_G(S(2m+1,2,1))$ is also determined.
\end{enumerate}
\end{lemma}
\begin{proof}
    We will only prove (1).
    Since
\[
m_1+m_2+m_3+m_4+t\leq 2m+4
\]
and $m_1,m_2,m_3,m_4,t\geq1$, at most one of
\[
m_1,m_2,m_3,m_4,t
\]
can be greater than or equal to $m+1$. Moreover, every one of these parameters is at most $2m$. 

Since $t \le 2m$, the value of $\delta_t(b)$ is determined by Lemma~\ref{lem:num_four_leaves}. Assume that $\delta_t(b) = 1$.

If $\sum_{i=1}^{2m} \gamma_i \ge 1$, after replacing every parameter greater than or equal to $2m+1$ by $2m+1$, the values of $a, c, k,$ and $\ell$ are uniquely determined.
Thus, we may assume that $\sum_{i=1}^{2m} \gamma_i = 0$. 

From $\alpha_1, \dots, \alpha_{2m}$, after replacing every parameter greater than or equal to $2m+1$ by $2m+1$, the multiset $\{a, c, k, \ell\}$, denoted by $\{x, y, z, w\}$, is determined. Let $x \ge y \ge z \ge w$. 
First, suppose that only $w$ is less than or equal to $m$. 
Then, it suffices to show that $G \cong \Pi(x, b, y, z, w)$ and $G' \cong \Pi(x, b, y, w, z)$ yield the same expression in Lemma~\ref{lem:num_four_leaves}. Since the two expressions are trivially identical if $U_{m_1}(x) = 0$ or $U_{m_1}(z) = 1$, we may assume that $U_{m_1}(x) = 1$ and $U_{m_1}(z) = 0$. 

Evaluating the formula in Lemma~\ref{lem:num_four_leaves} for $G$ yields
\[
U_{m_2}(w) + 2U_{m_3}(w) U_{m_1}(y) + U_{m_4}(w),
\]
and evaluating the same formula for $G'$ yields
\[
U_{m_2}(w) U_{m_1}(y) + 2U_{m_3}(w) + U_{m_4}(w) U_{m_1}(y).
\]
For these two values to differ, we must have $U_{m_1}(y) = 0$ and $U_{m_2}(w) + U_{m_4}(w) \neq 2U_{m_3}(w)$.

However, in this case, if $w < m_i$ for some $i \in \{1, 2, 3\}$, then $y + w \le 2m$, which contradicts $\sum_{i=1}^{2m} \gamma_i = 0$. Thus, if $U_{m_1}(y) = 0$, the condition $U_{m_2}(w) + U_{m_4}(w) \neq 2U_{m_3}(w)$ cannot be satisfied, and hence the expressions for $G$ and $G'$ coincide. 

The remaining cases, where the number of parameters among $x, y, z, w$ that are less than or equal to $m$ is $0$ or at least $2$, can be handled in an analogous manner.
\end{proof}

By Remark~\ref{rmk:a,c,k,l} and Lemmas~\ref{lem:A,B}, \ref{lem:Pi and T determined}, we have seen that the values of $\alpha_i$, $\beta_i$, and $\gamma_i$ play a key role in determining the structure of a $\Pi$-shape graph. 
We now show that the values of $\alpha_i$, $\beta_i$, and $\gamma_i$ needed for our argument can be recovered from the $2$-spectrum. 
We first show that these values can be determined from the numbers of $P_j$ and $T_j$ subgraphs.

\begin{lemma}\label{lem:recovery}
For each integer $m \geq 1$, suppose that the values $N_G(P_{j+3})$ and $N_G(T_{j+3})$ are given for all $j\leq m$.  
Then the quantities $\alpha_m$, $\beta_m$, and $\gamma_m$ are uniquely determined. 
\end{lemma}

\begin{proof}
We use an induction on $m$.

For $m=1$, we have
\[
N_G(P_4) = n + 1 - \alpha_1 + \beta_1, \qquad 
N_G(T_4) = 6 - \alpha_1 + 2\beta_1.
\]
Since $n$ is fixed, these two equations uniquely determine $\alpha_1$ and $\beta_1$.
Moreover, $\gamma_1 = 0$ trivially.

Assume that the statement holds for all integers less than $m$. 
Then $\alpha_i$, $\beta_i$, and $\gamma_i$ are known for all $1 \leq i \leq m-1$. 
Substituting these into the formulas of Lemmas~\ref{lem:num of path}, \ref{lem:num_tshape} and \ref{lem:A,B}, 
the expressions for $N_G(P_{m+3})$ and $N_G(T_{m+3})$ reduce to
\[
-\alpha_m + \beta_m + \gamma_m
\quad \text{and} \quad
-\alpha_m + 2\beta_m,
\]
respectively.

Taking the difference, we obtain $\beta_m - \gamma_m$. 
This determines $(\beta_m,\gamma_m)$ uniquely except in the cases where the difference coincides. 
There are two such possibilities, namely 
\[
(0,0) \text{ and } (1,1), \qquad \text{and} \qquad (1,2) \text{ and } (0,1).
\]

We first consider the pair $(0,0)$ and $(1,1)$. 
If both cases yield the same value of $-\alpha_m + \beta_m + \gamma_m$, then necessarily $\alpha'_m = \alpha_m + 2$, 
where $\alpha_m$ and $\alpha'_m$ correspond to the cases
$(\beta_m,\gamma_m)=(0,0), (1,1)$, respectively.
If $\alpha'_m\geq 3$, then $\gamma'_m=0$, and hence $\alpha_m=0$ and $\alpha'_m=2$. 
If $(\alpha'_m,\beta'_m,\gamma'_m)=(2,1,1)$, then the only possibility is
\[
\Pi(m,m,x,m,y), \qquad x+y=m,
\]
which has $4m$ edges.
On the other hand, if $(\alpha_m,\beta_m,\gamma_m)=(0,0,0)$, then the graph is of the form
\[
\Pi(a,b,c,x,y), \qquad a,b,c\geq m+1,
\]
and hence has at least $4m+3$ edges.
Therefore, the two cases have different numbers of edges and can be distinguished.

Next, consider the pair $(1,2)$ and $(0,1)$.
When $(\beta_m,\gamma_m)=(1,2)$, the total number of edges is $3m$, whereas when $(\beta_m,\gamma_m)=(0,1)$, it is at least $3m+2$.
This leads to an inconsistency with the fixed number of vertices of $G$, and hence cannot occur.

Therefore, all ambiguous cases are excluded, and $(\beta_m,\gamma_m)$ is uniquely determined. 
Consequently, $\alpha_m$ is also uniquely determined.
\end{proof}

We now show that the numbers of $P_j$ and $T_j$ subgraphs can be recovered from the $2$-spectrum for every $j$.

\begin{proposition}\label{prop:trace_even}
Let $G$ and $G'$ be $\Pi$-shape graphs. 
Assume that
\[
N_G(T_i)=N_{G'}(T_i) \quad \text{and} \quad N_G(P_i)=N_{G'}(P_i)
\]
for all $i=1,\dots,2m+3$. 
Suppose that
\[
\mathrm{tr}\big((A_G^{(2)})^{2m+4}\big)
=
\mathrm{tr}\big((A_{G'}^{(2)})^{2m+4}\big).
\]
Then
\[
N_G(P_{2m+4})=N_{G'}(P_{2m+4})
\quad \text{and} \quad
N_G(P_{2m+5})=N_{G'}(P_{2m+5}).
\]
\end{proposition}
\begin{proof}
Assume that 
\[
N_G(T_i)=N_{G'}(T_i), \quad N_G(P_i)=N_{G'}(P_i)
\]
for all $i=1,\dots,2m+3$. 
By Remark~\ref{rmk:a,c,k,l} and Lemma~\ref{lem:recovery}, for each of $a,c,k,\ell$, we can either determine its exact value if it is at most $m$, or determine that it is greater than $m$.
By Lemma~\ref{lem:Pi and T determined}, among the subgraphs contributing to 
\[
\mathrm{tr}\big((A^{(2)})^{2m+4}\big),
\]
all subgraph counts are uniquely determined except possibly 
\[
P_{2m+4}, \, P_{2m+5},\text{ and }S(m_1,m_2,m_3),
\]
where $2m \geq m_1 \geq m+1$ and $m \geq m_2 \geq m_3$.

We first show that the number of subgraphs of type $S(m_1,m_2,m_3)$ is uniquely determined.
If at most one of $a,c,k,\ell$ is at least $m+1$, then $N(S(m_1,m_2,m_3))$ is clearly determined by the known values of $\alpha_i$, $\beta_i$, and $\gamma_i$.
Thus, it suffices to consider the remaining cases.

\medskip
\noindent
\textbf{(1) Exactly two of $a,c,k,\ell$ are at least $m+1$.}

\smallskip
\noindent
\emph{(i) $a,c \geq m+1$.} 
Then $k,\ell\leq m$, so the multiset $\{k,\ell\}$ is determined. 
Since $a\ge c$, the values of $\alpha_{m+1},\dots,\alpha_{2m}$ allow us,
for each of $a$ and $c$, to either determine its exact value if it is at most $2m$, or determine that it is greater than $2m$.
While any parameter that is at least $2m+1$ may be replaced by $2m+1$ without changing the relevant unit step function values.
Thus, the values of $U_{m_1}(a)$ and $U_{m_1}(c)$ are determined.

If $k=\ell$, then $N(S(m_1,m_2,m_3))$ is uniquely determined. 
Hence, we may assume that $k\neq\ell$.
Thus, it suffices to compare $N(S(m_1,m_2,m_3))$ for
\[
G \cong \Pi(a,b,c,k,\ell)
\quad\text{and}\quad
G' \cong \Pi(a,b,c,\ell,k).
\]
Note that we use the same parameter symbols for $G$ and $G'$ to denote the corresponding truncated parameters, where every parameter at least $2m+1$ is replaced by $2m+1$. 
This does not affect any of the unit step function values relevant to the subgraphs under consideration.
Thus, the parameters denoted by the same symbol need not have the same original value in $G$ and $G'$. 
Whenever the original parameter values need to be distinguished, we use primes for the parameters of $G'$.

The values of $b$ and $b'$ can differ only if both exceed $2m$. 
In this case, however, the values of the unit step functions required to compute $N(S(m_1,m_2,m_3))$ are the same for $b$ and $b'$. 
Hence, we may assume that $b=b'$.

Since the values of $\gamma_1,\dots,\gamma_{2m}$ are the same for $G$ and
$G'$, we have $c+\ell\geq 2m+1$.

By Lemma~\ref{lem:num_tshape}, we obtain 
\[
\begin{aligned}
&N_G(S(m_1,m_2,m_3)) \neq N_{G'}(S(m_1,m_2,m_3)) \\
&\quad \iff
(U_{m_1}(a)-U_{m_1}(c)) \Big[
(U_{m_2}(k)-U_{m_2}(\ell))(1+U_{m_3}(b+\ell)-U_{m_3}(b)) \\
&\qquad\qquad\qquad +
(U_{m_3}(k)-U_{m_3}(\ell))(1+U_{m_2}(b+\ell)-U_{m_2}(b))
\Big] \\
&\qquad \neq 
(U_{m_1}(a+b)-U_{m_1}(b+c))
\big(U_{m_3}(k)+U_{m_2}(k)-U_{m_3}(\ell)-U_{m_2}(\ell)\big).
\end{aligned}
\]

If $c \geq m_1$ or $\ell \geq m_2$, then both sides vanish. 
Hence, the two values can differ only if $c < m_1$ and $\ell < m_2$. 
However, this implies
\[
c+\ell \leq m_1 + m_2 - 2 \leq 2m,
\]
which contradicts $c+\ell \geq 2m+1$. 
Hence the two values must coincide.

\smallskip
\noindent
\emph{(ii) $a,k \geq m+1$.} 
Then $c,\ell\leq m$, so both are already known. 
Since $a\ge c$, the values of $\alpha_{m+1},\dots,\alpha_{2m}$ determine, for each of $a$ and $c$, its exact value if it is at most $2m$, and otherwise determine that it is greater than $2m$.
If either $a$ or $k$ is at least $2m+1$, its exact value is irrelevant to $N(S(m_1,m_2,m_3))$, since $m_1\leq 2m$.
Thus, the values of $U_{m_1}(a)$ and $U_{m_1}(k)$ are uniquely determined,
and hence so is $N(S(m_1,m_2,m_3))$.

\medskip
\noindent
\textbf{(2) Exactly three of $a,c,k,\ell$ are at least $m+1$.}

Since $m_1\leq 2m$, we may replace every parameter that is at least $2m+1$ by $2m+1$ without changing any of the unit step function values appearing in $N(S(m_1,m_2,m_3))$. 
With this replacement, the values of $\alpha_{m+1},\dots,\alpha_{2m}$ determine the multiset of the three parameters among $a,c,k,\ell$ that are at least $m+1$, 
although they may not determine which three parameters among $a,c,k,\ell$ form this multiset.
Let this multiset be
\[
\{x,y,z\}, \qquad \text{with } \, x\geq y\geq z,
\]
and denote the remaining parameter by $w$.
Then, up to the symmetries of the $\Pi$-shape, there are three possible arrangements:
\[
\Pi(x,b,y,z,w),\qquad
\Pi(x,b,y,w,z),\qquad
\Pi(x,b,z,y,w).
\]
It is enough to compare the first two cases. 
Indeed, the two graphs
\[
\Pi(x,b,y,z,w)
\quad\text{and}\quad
\Pi(x,b,y,w,z)
\]
differ only by exchanging $x$ and $y$. 
Since the positions of $x$ and $z$ in $\Pi(x,b,y,z,w)$ can be exchanged by a symmetry of the $\Pi$-shape, 
the same comparison applies to any pair of the three parameters $x,y,z$. Hence, it suffices to consider
\[
G\cong\Pi(x,b,y,z,w),\qquad
G'\cong\Pi(x,b,y,w,z).
\]
Since the values of $\gamma_1,\dots,\gamma_{2m}$ are the same for $G$ and
$G'$, we have $y+w\geq 2m+1$.

By Lemma~\ref{lem:num_tshape},
\[
\begin{aligned}
&N_G(S(m_1,m_2,m_3)) \neq N_{G'}(S(m_1,m_2,m_3)) \\
&\quad \iff 
(U_{m_1}(x)-U_{m_1}(y)) \Big[
2 + U_{m_2}(b+w) - U_{m_2}(b) + U_{m_3}(b+w) - U_{m_3}(b) \\
&\qquad\qquad - U_{m_2}(w)(2 - U_{m_3}(b)) - U_{m_3}(w)(2 - U_{m_2}(b))
\Big] \\
&\qquad \neq 
(U_{m_1}(x+b)-U_{m_1}(b+y))(2 - U_{m_2}(w) - U_{m_3}(w)).
\end{aligned}
\]

If $U_{m_2}(w)=1$, then both sides are zero. 
If $U_{m_2}(w)=0$, then $y+w \geq 2m+1$, so both sides vanish. 
Hence the two values coincide.

\medskip
\noindent
\textbf{(3) All four of $a,c,k,\ell$ are at least $m+1$.}

In this case, all unit step function values appearing in the expression for $N(S(m_1,m_2,m_3))$, except possibly those involving $b$, are equal to $1$.
Since $U_{m_1}(a), U_{m_1}(c), U_{m_1}(k), U_{m_1}(\ell)$ appear symmetrically, 
the value of $N(S(m_1,m_2,m_3))$ is uniquely determined.

Thus, in all cases, the number of $S(m_1,m_2,m_3)$ is determined.

\medskip
Finally, we separate $N(P_{2m+4})$ and $N(P_{2m+5})$. 
After removing the known terms from the expression for $\mathrm{tr}((A^{(2)})^{2m+4})$, 
by Corollary~\ref{cor:first occurrence}, 
the remaining contribution is
\[
16^2 (m+2)^2 \, N(P_{2m+4}) + N(P_{2m+5}).
\]

By Lemmas~\ref{lem:num of path} and \ref{lem:A,B}, the undetermined part of $N(P_{2m+5})$ is
\[
-\alpha_{2m+2} + \beta_{2m+2} - \gamma_{2m+2},
\]
whose absolute value is strictly smaller than $16^2 (m+2)^2$. 
Therefore, the two values $N(P_{2m+4})$ and $N(P_{2m+5})$ are uniquely determined.
\end{proof}

\begin{proposition}\label{prop:trace_odd}
Let $G$ and $G'$ be $\Pi$-shape graphs.
Assume that
\[
N_G(T_i)=N_{G'}(T_i), \quad N_G(P_i)=N_{G'}(P_i)
\]
for all $i=1,\dots,2m+3$, and that
\[
N_G(P_{2m+4})=N_{G'}(P_{2m+4}), \quad 
N_G(P_{2m+5})=N_{G'}(P_{2m+5}).
\]
Suppose that
\[
\mathrm{tr}\big((A_G^{(2)})^{2m+5}\big)
=
\mathrm{tr}\big((A_{G'}^{(2)})^{2m+5}\big).
\]
Then
\[
N_G(T_{2m+4})=N_{G'}(T_{2m+4})
\quad \text{and} \quad
N_G(T_{2m+5})=N_{G'}(T_{2m+5}),
\]
except in the case where $m$ is odd and
\[
G\cong\Pi(a+c-2m-2,b,m+1,2m+2,m+1),\qquad
G'\cong\Pi(a,b,c,m+1,m+1)
\]
with $1\leq b\leq m+1$,and $a,c\geq2m+3$.

Furthermore, in this exceptional case, we have
\[
N_G(T_{2m+5})+1=N_{G'}(T_{2m+5}).
\]
\end{proposition}

\begin{proof}
As in the proof of Proposition~\ref{prop:trace_even}, among the subgraphs contributing to
$\mathrm{tr}\big((A^{(2)})^{2m+5}\big)$,
all subgraph counts are uniquely determined except possibly those of
\[T_{2m+4}, \quad T_{2m+5}, \quad
S(m+2,m+2,1), \quad S(m+2,m+1,1), \quad
S(m_1,m_2,m_3),\]
where $2m\geq m_1\geq m+1$ and $m\geq m_2\geq m_3$.
Note that $S(m+2,m+2,1)$ appears only when $m$ is odd.

We classify the remaining subgraphs into three types:
\begin{itemize}
\item Type 1: $S(m+2,m+2,1)$ (only when $m$ is odd),
\item Type 2: $S(m+2,m+1,1)$,
\item Type 3: $S(m_1,m_2,m_3)$ with $2m\geq m_1\geq m+1$ and $m\geq m_2\geq m_3$.
\end{itemize}

We analyze the contributions of these types to the trace and show that any possible discrepancies are incompatible with equality of the traces,
except in the exceptional case stated in the proposition.
The computations are analogous to those in Proposition~\ref{prop:trace_even}, and we omit repetitive details.

For Type 3, if $m_3\geq2$, the counts coincide by the same argument as before.
Thus, it suffices to consider $m_3=1$, which we further divide into
\begin{itemize}
\item Type 3-1: $S(m_1,m_2,1)$ with $m_1+m_2=2m+3$,
\item Type 3-2: $S(m_1,m_2,1)$ with $m_1+m_2=2m+4$ and $m_1,m_2$ are odd.
\end{itemize}

We proceed according to the number of elements among $a,c,k,\ell$ that are at least $m+1$.
\medskip
\noindent
\textbf{(1) Exactly two of $a,c,k,\ell$ are at least $m+1$.}
As in the previous proof, we may assume
\[G\cong\Pi(a,b,c,k,\ell),\qquad
G'\cong\Pi(a,b,c,\ell,k),\]
with $a>c\geq m+1$ and $m\geq k>\ell$.

Recall that we use the same parameter symbols for $G$ and $G'$ to denote the corresponding truncated parameters, where every parameter at least $2m+1$ is replaced by $2m+1$. 
Whenever the original parameter values need to be distinguished, we use primes for the parameters of $G'$.

Applying Lemma~\ref{lem:num_tshape}, one checks that the counts of Type 1 and Type 2 always coincide.
For Type 3, the possible discrepancies are
$$
\begin{array}{c|c|c}
\text{Type} & \text{Condition} & N_G-N_{G'} \\ \hline
3\text{-}1
& c=m_1-1,\ \ell=m_2-1
& 1 \\ \hline
3\text{-}2
& c=m_1-1,\ \ell=m_2-1
& 1 \\
3\text{-}2
& c=m_1-1,\ \ell=m_2-2,\ b=1
& -1 \\
3\text{-}2
& c=m_1-1,\ \ell=m_2-2,\ b\geq2
& 1 \\
3\text{-}2
& c=m_1-2,\ \ell=m_2-1,\ a=m_1-1
& -1
\end{array}
$$
and the counts coincide in all other cases.

Since
\[
\mathrm{tr}\big((A_G^{(2)})^{2m+5}\big)
=
\mathrm{tr}\big((A_{G'}^{(2)})^{2m+5}\big),
\]
Corollary~\ref{cor:first occurrence} implies that
\[
\begin{aligned}
&\big\{N(T_{2m+5})+N(S_2(m_1,m_2,1))\big\}\\
&\qquad
+16^2\big\{(m+2)^2N(T_{2m+4})
+(m+2)(m_i/2+1)N(S_1(m_1,m_2,1))\big\}
\end{aligned}
\]
has the same value for $G$ and $G'$, 
where $S_1(m_1,m_2,1)$ and $S_2(m_1,m_2,1)$ are of Type 3-1 and Type 3-2, respectively, 
and $m_i$ denotes the even one among $m_1$ and $m_2$ in Type 3-1.

By assumption, $\alpha_j$, $\beta_j$, and $B_j$ are the same for $G$ and $G'$ for $j=1,\dots,2m$. 
Hence, by Lemma~\ref{lem:num_tshape}, the difference between the values of $N(T_{2m+5})$ for $G$ and $G'$ is at most $8$. 
On the other hand, from the table above, the difference between the values of $N(S_2(m_1,m_2,m_3))$ is at most $1$. 
Therefore, the difference between
\[
N(T_{2m+5})+N(S_2(m_1,m_2,m_3))
\]
is at most $9$, whereas the second term is multiplied by $16^2$. It follows that the two quantities
\[
N(T_{2m+5})+N(S_2(m_1,m_2,m_3)) \quad \text{and }\,
(m+2)^2N(T_{2m+4})
+(m+2)(m_i/2+1)N(S_1(m_1,m_2,m_3))
\]
must have the same value for $G$ and $G'$.

First, suppose that the numbers of Type 3-1 subgraphs are different. 
Then, from
\[
(m+2)^2N(T_{2m+4})
+(m+2)(m_i/2+1)N(S_1(m_1,m_2,m_3)),
\]
we obtain
\[
m_i/2+1
=
(m+2)\bigl(N_{G'}(T_{2m+4})-N_G(T_{2m+4})\bigr).
\]
However, $m+2$ cannot divide $m_i/2+1$. 
Thus, the above equality is impossible. Hence, such a pair of graphs cannot satisfy
$
\mathrm{tr}\big((A_G^{(2)})^{2m+5}\big)
=
\mathrm{tr}\big((A_{G'}^{(2)})^{2m+5}\big).
$

Next, suppose that the numbers of Type 3-2 subgraphs are different. 

Consider first the case
$c=m_1-1,\qquad \ell=m_2-1.$
The value of $N(T_{2m+4})$ determines $\alpha_{2m+1}$ and $\beta_{2m+1}$, while the equality of
\[
N(T_{2m+5})+N(S_2(m_1,m_2,m_3))
\]
gives
\[
1=N_{G'}(T_{2m+5})-N_G(T_{2m+5}).
\]
By Lemma~\ref{lem:num_tshape}, the latter equality can hold only when
\[
b\geq 2m+3,\qquad b'=2m+2,\qquad
\alpha_{2m+2}=0,\qquad
\alpha'_{2m+2}=1.
\]
But in this case $G$ and $G'$ have different numbers of edges, a contradiction. 
The remaining three cases can be ruled out in the same way.

Therefore, whenever the number of subgraphs of at least one of Types 1, 2, and 3 differs between $G$ and $G'$, we have
$
\mathrm{tr}\big((A_G^{(2)})^{2m+5}\big)
\neq
\mathrm{tr}\big((A_{G'}^{(2)})^{2m+5}\big).
$
In all remaining cases, the same argument as in the proof of Proposition~\ref{prop:trace_even} yields
\[
N_G(T_{2m+4})=N_{G'}(T_{2m+4})
\quad\text{and}\quad
N_G(T_{2m+5})=N_{G'}(T_{2m+5}).
\]

\medskip
\noindent
\textbf{(2) Exactly three of $a,c,k,\ell$ are at least $m+1$.}
Assume $x,y,z\geq m+1$ and $w\leq m$, and compare
\[G\cong\Pi(x,b,y,z,w),\qquad
G'\cong\Pi(x,b,y,w,z).\]
The possible discrepancies for Types 1, 2, and 3 are

$$
\begin{array}{c|c|c}
\text{Type} & \text{Condition} & N_G-N_{G'} \\ \hline
1
& b=1,\ x\geq m+2,\ y=z=m+1,\ w=m
& -1 \\
1
& b\geq2,\ x,z\geq m+2,\ y=m+1,\ w=m
& 1 \\ \hline
2
& b=1,\ x\geq m+2,\ y=z=m+1,\ w=m
& -2 \\
2
& 2\leq b\leq m+1,\ x\geq m+2,\ y=z=m+1,\ w=m
& -3 \\
2
& b\geq m+2,\ x\geq m+2,\ y=z=m+1,\ w=m
& -2 \\
2
& x,z\geq m+2,\ y=m+1,\ w=m
& 1 \\ \hline
3\text{-}1
& z=m_1-1,\ w=m_2-1
& -1 \\ \hline
3\text{-}2
& z=m_1-1,\ w=m_2-1
& -1 \\
3\text{-}2
& z=m_1-1,\ w=m_2-2,\ b\geq2
& 1 \\
3\text{-}2
& x=m_1-1,\ z=m_1-2,\ w=m_2-1,\ b=1
& -1 \\
3\text{-}2
& x\geq m_1,\ z=m_1-2,\ w=m_2-1,\ b\geq2
& 1
\end{array}
$$
and the counts coincide in all other cases.
As in Case (1), it follows that the two quantities
\[
\begin{aligned}
&N(T_{2m+5})+N(S_2(m_1,m_2,m_3))
+N(S(m+2,m+2,1))
\end{aligned}
\]
and
\[
\begin{aligned}
&(m+2)N(T_{2m+4})
+(m_i/2+1)N(S_1(m_1,m_2,m_3))+(\left\lfloor\frac{m+2}{2}\right\rfloor+1)N(S(m+2,m+1,1))
\end{aligned}
\]
must have the same value for $G$ and $G'$.

For each of the cases in which the number of subgraphs of at least one of the three types differs between $G$ and $G'$, 
the calculation is identical to that in Case (1), and leads to a contradiction. 
Therefore, whenever the number of subgraphs of at least one of Types 1, 2, and 3 differs between $G$ and $G'$, we have

$$
\mathrm{tr}\big((A_G^{(2)})^{2m+5}\big)
\neq
\mathrm{tr}\big((A_{G'}^{(2)})^{2m+5}\big).
$$

\medskip
\noindent
\textbf{(3) All of $a,c,k,\ell$ are at least $m+1$.}
Again, compare
$$
G\cong\Pi(x,b,z,y,w),\qquad
G'\cong\Pi(x,b,y,w,z),
$$
with $x\ge y \ge z \ge w \ge m+1$.

The only possible discrepancies are
$$
\begin{array}{c|c|c}
\text{Type} & \text{Condition} & N_G-N_{G'} \\ \hline
1
& m\text{ odd},\ x,y\geq m+2,\ z=w=m+1
& 1 \\
2
& x,y\geq m+2,\ z=w=m+1,\ b\geq m+2
& -2
\end{array}
$$
and the counts coincide in all other cases.

The Type 2 discrepancy leads to a contradiction as before.
Thus we may assume that $1\leq b\leq m+1$.

For a Type 1 discrepancy, as in Case (1), the equality
\[
\mathrm{tr}\big((A_G^{(2)})^{2m+5}\big)
=
\mathrm{tr}\big((A_{G'}^{(2)})^{2m+5}\big)
\]
and Corollary~\ref{cor:first occurrence} imply that
\[
N(T_{2m+5})+N(S(m+2,m+2,1))
\]
has the same value for $G$ and $G'$. 

Hence,
\[
1=N_{G'}(T_{2m+5})-N_G(T_{2m+5}).
\]

By Lemma~\ref{lem:num_tshape}, this gives
\[
1=\alpha_{2m+2}-\alpha'_{2m+2}.
\]

This can occur only when
\[
(\alpha_{2m+2},\alpha'_{2m+2})=(1,0).
\]

Consequently, the only exceptional case is
\[
G\cong\Pi(a+c-2m-2,b,m+1,2m+2,m+1),\qquad
G'\cong\Pi(a,b,c,m+1,m+1)
\]
where
\[
1\leq b\leq m+1,\qquad a,c\geq2m+3,
\]
and $m$ is odd. 
Excluding this exceptional case, we have
\[N_G(T_{2m+4})=N_{G'}(T_{2m+4})
\quad\text{and}\quad
N_G(T_{2m+5})=N_{G'}(T_{2m+5}).
\]
\end{proof}

\begin{proposition}\label{prop:final_distinguish}
Suppose that $m$ is odd, $1\leq b\leq m+1$, and $a,c\geq 2m+3$.
Let 
\[
G\cong\Pi(a+c-2m-2,b,m+1,2m+2,m+1),\qquad
G'\cong\Pi(a,b,c,m+1,m+1).
\]
Then
\[
\mathrm{tr}\big((A_G^{(2)})^{2m+6}\big)
\neq 
\mathrm{tr}\big((A_{G'}^{(2)})^{2m+6}\big).
\]
Thus, $G$ and $G'$ are not $2$-cospectral.
\end{proposition}
\begin{proof}
Assume, for a contradiction, that
\[
\mathrm{tr}\big((A_G^{(2)})^{2m+6}\big)
=
\mathrm{tr}\big((A_{G'}^{(2)})^{2m+6}\big).
\]

By Corollary~\ref{cor:first occurrence}, the subgraphs that first appear and need to be considered are of the following types:
\begin{enumerate}
    \item $S(2m+2,2,1)$;
    \item $S(m+2,m+1,2)$;
    \item $S(m_1,m_2,m_3)$, where
    \[
        m_1+m_2+m_3=2m+5
        \quad\text{or}\quad
        m_1+m_2+m_3=2m+6,
    \]
    with
    \[
        2m+1\geq m_1\geq m+2,
        \qquad
        m+1\geq m_2 \geq m_3.
    \]
    \item $P_{2m+6}$ and $P_{2m+7}$.
\end{enumerate}

We first compare the numbers of the three types of $T$-shape
subgraphs using Lemma~\ref{lem:num_tshape}.
For $S(2m+2,2,1)$, we have
\begin{align*}
N_G(S(2m+2,2,1))
&=
1+2(2-U_2(b))+1+4U_{m+1}(b)+4\\
&=
14-2U_2(b)+4U_{m+1}(b),
\end{align*}
whereas
\begin{align*}
N_{G'}(S(2m+2,2,1))
&=
2\big\{1+2-U_2(b)+2(1+U_{m+1}(b))\big\}\\
&=
10-2U_2(b)+4U_{m+1}(b).
\end{align*}
Thus,
\[
N_G(S(2m+2,2,1))
=
N_{G'}(S(2m+2,2,1))+4.
\]

For $S(m+2,m+1,2)$, we have
\[
N_G(S(m+2,m+1,2))
=
N_{G'}(S(m+2,m+1,2))
=
16-2U_2(b)-2U_{m+1}(b).
\]

Finally, for $S(m_1,m_2,m_3)$ satisfying the conditions in (3), we have
\[
N_G(S(m_1,m_2,m_3))
=
N_{G'}(S(m_1,m_2,m_3))
=
10-2U_3(b)-2U_2(b)
+4U_{m_1}(b+m+1).
\]

As in the case of $\mathrm{tr}\big((A_G^{(2)})^{2m+5}\big)$, the equality
\[
\mathrm{tr}\big((A_G^{(2)})^{2m+6}\big)
=
\mathrm{tr}\big((A_{G'}^{(2)})^{2m+6}\big)
\]
allows us to distinguish the contributions from subgraphs with $2m+5$ and $2m+6$ edges. 
Hence, among the subgraphs with $2m+5$ edges, the only ones whose numbers of subgraphs differ between $G$ and $G'$ are
\[
S(2m+2,2,1),\qquad
S(m+2,m+2,1),\qquad
T_{2m+5},\qquad
P_{2m+6}.
\]

Moreover, by the corresponding formulas for the numbers of these
subgraphs,
\[
N_G(S(m+2,m+2,1))
=
N_{G'}(S(m+2,m+2,1))+1
\]
and
\[
N_G(T_{2m+5})
=
N_{G'}(T_{2m+5})-1.
\]

Taking the coefficients of these four subgraphs in the expansion of
$\mathrm{tr}((A^{(2)})^{2m+6})$ into account, the assumed equality of
the two traces gives
\begin{align*}
&2(m+2)N_G(S(2m+2,2,1))
+(m+3)^2N_G(P_{2m+6})\\
&\quad+
\left(
\frac{(m+3)(m+2)}{2}
+
\frac{(m+5)(m+3)}{4}
\right)
N_G(S(m+2,m+2,1))\\
&\quad+
\big((m+3)(m+2)+1\big)N_G(T_{2m+5})\\
={}&
2(m+2)N_{G'}(S(2m+2,2,1))
+(m+3)^2N_{G'}(P_{2m+6})\\
&\quad+
\left(
\frac{(m+3)(m+2)}{2}
+
\frac{(m+5)(m+3)}{4}
\right)
N_{G'}(S(m+2,m+2,1))\\
&\quad+
\big((m+3)(m+2)+1\big)N_{G'}(T_{2m+5}).
\end{align*}
Hence,
\begin{align*}
&8(m+2)
+\frac{(m+3)(m+2)}{2}
+\frac{(m+5)(m+3)}{4}
-\big((m+3)(m+2)+1\big)\\
&\quad=(m+3)^2
\big(
N_{G'}(P_{2m+6})-N_G(P_{2m+6})
\big).
\end{align*}

Thus, by Lemma~\ref{lem:num of path}, we obtain
\begin{equation}\label{equation}
\frac{-m^2+30m+63}{4}
=
(m+3)^2
\big(
\alpha_{2m+3}-\alpha'_{2m+3}
\big).    
\end{equation}

There are only three possible pairs:
\[
(\alpha_{2m+3},\alpha'_{2m+3})
=
(0,0),\qquad
(0,1),\qquad
(0,2).
\]

For the three pairs above, equation~\eqref{equation} reduces, respectively, to
\[ -m^2+30m+63=0,\qquad
3m^2+54m+99=0,\qquad
7m^2+78m+135=0. \]
None of these equations has a positive integer solution.

Therefore,
\[
\mathrm{tr}\big((A_G^{(2)})^{2m+6}\big)
\neq
\mathrm{tr}\big((A_{G'}^{(2)})^{2m+6}\big).
\]
Hence, $G$ and $G'$ are not $2$-cospectral.
\end{proof}

The following corollary is an immediate consequence of Propositions \ref{prop:trace_even}, \ref{prop:trace_odd}, and \ref{prop:final_distinguish}.

\begin{corollary} \label{coro:P,T}
Let $G$ and $G'$ be $\Pi$-shape graphs on $n$ vertices. Then $G$ and $G'$ are $2$-cospectral if and only if
\[
N_G(P_m)=N_{G'}(P_m)
\quad\text{and}\quad
N_G(T_m)=N_{G'}(T_m)
\]
for all $m$.
\end{corollary}

\begin{theorem}
Every $\Pi$-shape graph is determined by its $2$-spectrum among the $C_3$-free connected graphs.
\end{theorem}
\begin{proof}
Suppose that $G$ and $G'$ are $2$-cospectral $C_3$-free graphs, where $G$ is a $\Pi$-shape graph.
By Lemmas~\ref{lem:determine H of tree} and~\ref{lem:determine H of C3 free}, $G'$ is also a $\Pi$-shape graph.
By Remark~\ref{rmk:a,c,k,l}, Lemma~\ref{lem:recovery},
and Corollary~\ref{coro:P,T}, the equality of the $2$-spectra of $G$ and $G'$ implies that the parameters $a,b,c,k,\ell$ of $G$ and $G'$ coincide.
Thus every $\Pi$-shape graph is determined by its $2$-spectrum among the $C_3$-free graphs.
\end{proof}

\bigskip

\appendix
\section{SageMath Code for Verification}
\label{app:sage_code}

The following SageMath script was used to verify that all graphs with at most nine vertices are determined by their $2$-spectra.
For each $n \le 9$, all non-isomorphic graphs of order $n$ were generated using \texttt{nauty\_geng}, and the characteristic polynomials of their $2$-adjacency matrices were compared.

\begin{lstlisting}
    from sage.graphs.graph_generators import graphs
from sage.matrix.constructor import diagonal_matrix

def second_adjacency_charpoly(G):
    A = G.adjacency_matrix()
    D = diagonal_matrix(G.degree())
    M = 16*A - A**2 + D
    return M.charpoly()


# Search for pairs of non-isomorphic graphs having
# the same 2-adjacency characteristic polynomial.
for n in range(1, 10):

    spectra = {}
    collision_found = False
    total_graphs = 0

    # nauty_geng generates one representative
    # from each isomorphism class of graphs of order n.
    for G in graphs.nauty_geng(str(n)):

        total_graphs += 1
        cp = second_adjacency_charpoly(G)

        if cp in spectra:
            collision_found = True
            break

        spectra[cp] = True

    status = (
        "COLLISION FOUND"
        if collision_found
        else "ALL GRAPHS ARE DISTINGUISHED"
    )

    print(
        f"n = {n} "
        f"({total_graphs:6d} graphs checked): "
        f"{status}"
    )
\end{lstlisting}

\bigskip

\section*{Statements and Declarations}

\subsection*{Funding}
The authors declare that no funds, grants, or other support were received during the preparation of this manuscript.

\subsection*{Competing Interests}
The authors have no relevant financial or non-financial interests to disclose.

\subsection*{Data Availability}
Not applicable.

\subsection*{Declaration of Generative AI Use}
ChatGPT (OpenAI, GPT-5.6) was used to assist with correcting English grammar and refining wording and expressions for clarity and readability. The authors reviewed and verified all suggested changes and take full responsibility for the content of the manuscript.

\end{document}